\documentclass[12pt]{amsart}
\usepackage{amsmath,amssymb,amsthm,amscd,amsopn,amsfonts,amsxtra} 
\usepackage{latexsym,array,verbatim,url}
\usepackage[mathscr]{eucal}
\usepackage{graphicx,psfrag,epsfig,color}

\usepackage[small,nohug,heads=littlevee]{diagrams} 
\diagramstyle[labelstyle=\scriptstyle]

\newtheorem{theorem}[subsection]{Theorem} 
\newtheorem{lemma}[subsection]{Lemma}

\newtheorem{assumption}[subsection]{Assumption}

\newtheorem{definition}[subsection]{Definition}

\def\ga{\gamma}

\def\veps{\varepsilon}
\def\vphi{\varphi}
\def\De{\Delta}
\def\Ga{\Gamma}

\def\supp{\mathrm{supp}}

\def\cL{\mathcal{L}}
\def\inv{^{-1}}
\def\iy{\infty}
\def\pa{\partial} 

\newcommand{\cal}{\mathcal}

\newcommand{\Z}{\mathbb{Z}}
\newcommand{\R}{\mathbb{R}}
\newcommand{\C}{\mathbb{C}}
\newcommand{\N}{\mathbb{N}}

\def\sideremark#1{\ifvmode\leavevmode\fi\vadjust{\vbox to0pt{\vss
 \hbox to 0pt{\hskip\hsize\hskip1em
\vbox{\hsize2cm\tiny\raggedright\pretolerance10000
 \noindent #1\hfill}\hss}\vbox to8pt{\vfil}\vss}}}%

\begin{document}
\title[Interpolation Theorems for Self-adjoint  operators]
{Interpolation Theorems for Self-adjoint  operators}
\author{Shijun Zheng}
\address{Department of Mathematical Sciences\\
Georgia Southern University\\
Statesboro, GA 30460-8093}
\email{szheng@georgiasouthern.edu}
\urladdr{http://math.georgiasouthern.edu/\symbol{126}{szheng}}
\date{\today}

\begin{abstract}  
We prove  a complex and a real interpolation \mbox{theorems} on Besov spaces and Triebel-Lizorkin spaces 
associated with a selfadjoint operator $\cL$, without assuming the gradient estimate for its spectral kernel.
The result applies to the cases where $\cL$ is a uniformly elliptic operator or a Schr\"odinger operator
with electro-magnetic potential. 
\end{abstract}

\maketitle 

\section{Introduction and main result}\label{S1} 
Interpolation of 
function spaces has  
played an important role in classical Fourier analysis and PDEs \cite{BS93,CuS01,
GV95,JN94,OZ08}. 
Let $\cL$ be a selfadjoint  operator in $L^2(\R^n)$. Then, for a Borel measurable function $\phi$: $\R\to \C$, 
we define $\phi(\cL)$ using functional calculus. 
 In \cite{JN94,D97,BZ,OZ06}  several authors introduced and studied the 
Besov spaces and Triebel-Lizorkin spaces associated with Schr\"odinger operators. In this note we 
present an 
interpolation 
result on these spaces for 
$\cL$.  

Let $\{\varphi_j\}_{j=0}^\iy\subset C_0^\infty ({\R}) $ be a dyadic system satisfying 
(i)  $\supp \,\vphi_0 \subset \{ x: |x|\le 1\}$,
 $\supp\; \varphi_j
\subset \{ x: 2^{j-2}\le |x|\le 2^{j}\} $, $j\ge 1$, 
(ii) $|\vphi_j^{(k)}(x)|\le c_k 2^{-kj}$ for all $j$, $k\in \N_0=\N\cup\{0\}$,
(iii) $\displaystyle{
\sum_{j=0}^\infty |\varphi_j(x)| \approx 1, \quad \forall x.}$
Let $\alpha\in \R$,  $1\le p,q\le \infty$. 
The inhomogeneous \emph{Besov space} \emph{associated with} $\cL$, denoted by ${B}_p^{\alpha,q}(\cL)$, %
is defined to be the completion of $\mathcal{S}(\R^n)$, the Schwartz class,  with respect to  
the norm 
$$ \Vert f\Vert_{{B}_p^{\alpha,q}(\cL) } = 
\big(\sum_{j=0}^{\infty}
 2^{j\alpha q} \Vert \vphi_j(\cL)f \Vert_{L^p}^q \big)^{1/q}\, .
$$
Similarly, the inhomogeneous {\em Triebel-Lizorkin space associated with $\cL$}, denoted by ${F}_p^{\alpha,q}(\cL)$, 
is defined by the norm
\begin{equation*}\label{eq:F-norm}
\Vert f\Vert_{{F}_p^{\alpha,q}(\cL)} 
= 
\Vert \big(\sum_{j=0}^{\infty} 2^{j\alpha q} \vert \vphi_j(\cL)f \vert^q\big)^{1/q}\Vert_{L^p} \, .
\end{equation*} 

The following assumption on the kernel of $\phi_j(\cL)$ is fundamental in the study of 
function space theory. 
 Let $\phi(\cL)(x,y)$ denote the integral 
kernel of $\phi(\cL)$. 
\begin{assumption}\label{a:phi-dec-Om}  
Let  $\phi_j\in C_0^\infty(\R)$ satisfy conditions (i), (ii) above.  
Assume that 
  there exist some $\veps>0$ 
and a constant 
 $c_n>0$ such that for all $j$
\begin{equation}\label{e:ker-phi-dec}
\vert  \phi_j(\cL)(x,y)\vert 
\le c_n\frac{2^{nj/2}}{(1+2^{j/2}|x-y|)^{n+\veps}}\,. 
\end{equation}
\end{assumption}
This is the same condition assumed in \cite{Z06a,OZ08} except that we drop the 
gradient estimate condition on the kernel. 
This is the case when $\cL$ is a Schr\"odinger operator $-\De+V$, 
$V\ge 0$ belonging to $L^1_{loc}(\R^n)$ \cite{He90a,ou06} or $\cL$ is a uniformly elliptic operator in $L^2(\R^n)$
\cite[Theorem 3.4.10]{Da89}. 

In what follows, $[A,B]_\theta$ denotes the usual complex interpolation between two Banach spaces; 
$(A,B)_{\theta,r}$ the real interpolation, see \mbox{Section \ref{S2}.} The notion 
 $T$: $X\to Y$ means that the \mbox{linear} operator $T$ is bounded from $X$ to $Y$.

\begin{theorem}[complex interpolation]\label{th:F-B-com-int}
Suppose that $\cL$ is a selfadjoint operator satisfying Assumption \ref{a:phi-dec-Om}. 
Let $0<\theta<1$,  $s=(1-\theta)s_0+\theta s_1 $, $s_0,s_1\in\R$ 
and  \begin{equation*}
\frac{1}{p}=\frac{1-\theta}{p_0}+ \frac{\theta}{ p_1},\quad 
\frac{1}{q}=\frac{1-\theta}{q_0}+ \frac{\theta}{ q_1}.  \end{equation*}
(a) If  $ 1< p_i<\infty$,  $1< q_i< \infty$, $i=0, 1$, then 
$$[{F}_{p_0}^{s_0,q_0}(\cL),{F}_{p_1}^{s_1,q_1}(\cL)]_\theta={F}_{p}^{s,q}(\cL). 
$$
(b) If $ 1\le p_i\le\iy$, $1\le q_i\le\infty$, $i=0,1$, 
 then  $$[{B}_{p_0}^{s_0,q_0}(\cL),{B}_{p_1}^{s_1,q_1}(\cL)]_\theta
={B}_{p}^{s,q}(\cL). 
$$
c) If $T: {F}_{p_0}^{s_0,q_0}(\cL)\to {F}_{\bar{p}_0}^{\bar{s}_0,\bar{q}_0}(\cL)$ 
and $T: {F}_{p_1}^{s_1,q_1}(\cL)\to {F}_{\bar{p}_1}^{\bar{s}_1,\bar{q}_1}(\cL)$, 
then $T: {F}_{p}^{s,q}(\cL) \to F_{\bar{p}}^{\bar{s},\bar{q}}(\cL) $,
where $\bar{s},\bar{p},\bar{q}$ and $\bar{s}_i,\bar{p}_i,\bar{q}_i$,  
satisfy the same relations as those for $s,p,q$ and $s_i,p_i,q_i$, $1<p_i,q_i<\iy$.
 Similar statement holds for $B^{s,q}_p(\cL)$.
\end{theorem}

 Complex  interpolation method
 originally was due to Calder\'on \cite{Ca63} and Lions and Peetre \cite{LP64}; see also \cite{Gr66,Taib}. 
The classical interpolation \mbox{theory} for Besov and Triebel-Lizorkin spaces on $\R^n$ has been given systematic treatments 
in \cite{Pe76}, \cite{BL76}, and \cite{Tr78,Tr83}. 
 There are  interesting discussions on interpolation theory  
in \cite{Pe76} and \cite{Tr83, Tr78, SchmTr87} for generalized Besov spaces associated with differential operators, which 
requires certain Riesz summability for 
$\cL$ that seems a nontrivial condition to verify. Nevertheless, we would like to mention that the Riesz summability, the spectral multiplier
theorem and the decay estimate in (\ref{e:ker-phi-dec}) are actually intimately related \cite{DOS02,OZ08}.

The real interpolation result for $B^{\alpha,q}_p(\R^n)$, $F^{\alpha,q}_p(\R^n)$ can be found in 
\cite{Pe76} and \cite{Tr83,Tr78}. Following the proof as in the classical case, 
but applying the estimate in 
(\ref{e:ker-phi-dec}) in stead of  spectral multiplier result, we obtain   
\begin{theorem}[real interpolation]\label{th:F-B-real-int}
Suppose that $\cL$ 
 satisfies Assumption \ref{a:phi-dec-Om}. Let 
\mbox{$0<\theta<1$,} $1\le r\le \infty$, 
  $s=(1-\theta)s_0+\theta s_1 $,  $s_0\neq s_1$. \\
(a) If $1\le p<\iy, 1\le q_1,q_2 \le\infty$, then
$$(F_{p}^{s_0,q_0}(\cL),F_{p}^{s_1,q_1}(\cL) )_{\theta,r}=B_{p}^{s,r}(\cL). 
$$
(b) If $ 1\le p,  q_1,q_2\le\infty$, 
then  $$(B_{p}^{s_0,q_0}(\cL),B_{p}^{s_1,q_1}(\cL) )_{\theta,r}
=B_{p}^{s,r}(\cL). 
$$ 
\end{theorem}







 The homogeneous spaces $\dot{B}_p^{\alpha,q}(\cL)$ and $\dot{F}_p^{\alpha,q}(\cL)$ can be defined using $\{\varphi_j\}_{j=-\iy}^\iy$ in (i) to (iii), 
instead of  $\{\vphi_j\}_{j=0}^\iy$. Then the analogous results of Theorem \ref{th:F-B-com-int} and Theorem \ref{th:F-B-real-int} hold. 

\section{interpolation for $\cL$}\label{S2} 


Theorem \ref{th:F-B-com-int} and Theorem \ref{th:F-B-real-int} are part of 
the abstract interpolation theory for $\cL$. 
In this section we present  the outline of their proofs. 
It was mentioned   in \cite{SchmTr87} that the interpolation associated with $\cL$  
is a ``subtle and difficult" subject, which normally relies on the very property of $\cL$.

\subsection{Complex interpolation} 
The proof of Theorem \ref{th:F-B-com-int} is similar to that
 given in \cite{Tr78} in the Fourier case. 
The insight is that the three line theorem (involved in Riesz-Thorin or Calder\'on's constructive proof
for $L^p$ spaces) 
 reflects the fact that the value of an 
 analytic function 
 in the interior of a domain is 
determined by its boundary values. 
\begin{definition}
Let $(A_0,A_1)$ be an interpolation couple , i.e.,
$A_0,A_1$ are (complex) Banach spaces, linearly and continuously 
embedded in a Hausdorff space $\cal{H}$. 
The space $A_0\cap A_1$ is endowed with the norm 
$\Vert a\Vert_{A_0\cap A_1}= 
\max\{\Vert a\Vert_{A_j},j=0,1\}$. The space $A:=A_0+A_1$ 
is endowed with the norm 
\begin{equation*}
\Vert a\Vert_{A}=\inf
\{\Vert a_0\Vert_{A_0}+\Vert a_1\Vert_{A_1}: a_0\in A_0,a_1\in A_1\} .
\end{equation*}


Let $S=\{z\in\C:0\le\Re z\le 1\}$ and $\bar{S}$ its closure. 
Denote $F$ the class of all $A$
-valued functions $f(z)$ on 
$\bar{S}$ such that $z\mapsto f(z)\in A$ is analytic in $S$ and continuous on $\bar{S}$,
satisfying
\begin{enumerate}
\item[(i)]  \[
\sup_{z\in \bar{S}} \Vert f(z)\Vert_A \;\text{is finite}. \]
\item[(ii)] The mapping $t\mapsto f(j+it)\in A_j$
 are continuous from $\R$ to $A_j$, $j=0,1$.\end{enumerate}
Then $F$ is a Banach space with the norm
\[
\Vert f\Vert_F= \max_j\{\sup_t  \Vert f(j+it)\Vert_{ A_j}\}.
\]
For $0<\theta<1$ we define the interpolation space $[A_0,A_1]_\theta$
as
\[ [A_0,A_1]_\theta:=\{a\in A: \exists f\in F
\;\text{with}\; f(\theta)=a\}  .\]
Then $[A_0,A_1]_\theta$ is a Banach space equipped with the norm
\[ \Vert a\Vert_\theta :=\inf\{ \Vert f\Vert_F: f\in F\, and\, f(\theta)=a\}.\]
\end{definition}








\subsection{Outline of the proof of Theorem \ref{th:F-B-com-int}} 
Let $\{\phi_j\},\{\psi_j\}$ satisfy the conditions in (i)-(iii) and
 $\sum_j \psi_j(x)\phi_j(x)=1$.
Define the operators $S: 
f\mapsto \{\phi_j(\cL)f\}$,  and $R:g\mapsto \sum_j\psi_j(\cL)g$.
The proof for part (a) follows from  the commutative diagram
$$\begin{CD}
F^{s,q}_p(\cL)  @>S >>L^{p}(\ell^q)\\ 
@VIdVV             @VVIdV\\    
F^{s,q}_p(\cL)  @<R<<L^p(\ell^q) 
\end{CD}
$$ 
and Lemma \ref{lem:vector-Lp-int} and
Lemma \ref{lem:vector-lq-int},  
which are  interpolation results for Banach space valued $L^p$ and $\ell^q$ spaces \cite{Tr78}.


\begin{lemma} \label{lem:vector-Lp-int}
Let $0<\theta<1$, $1\le p_0,p_1<\infty$ and $p\inv=(1-\theta)p_0\inv+\theta p_1\inv$.  Let $A_0,A_1$ be Banach spaces. Then 
\begin{equation}\label{eq:vector-Lp-int}
[L^{p_0}(A_0),L^{p_1}{(A_1)}]_\theta
=L^{p}([A_0,A_1]_\theta).   
\end{equation}
If $p_1=\infty$, 
then (\ref{eq:vector-Lp-int}) holds with 
$L^{p_1}(A_1)$ replaced by ${L_0^\infty}(A_1)$, the completion of simple $A_1$-valued functions with the esssup norm.
\end{lemma}

As in \cite{Tr78}, denote $\ell^q(A_j)$ the space of functions consisting of $a=\{a_j\}$, $a_j\in A_j$ ($A_j$ being Banach spaces)
equipped with the norm\\ 
$\displaystyle \Vert a\Vert_{\ell^q(A_j)}=\left( \sum_j\Vert a_j\Vert^q_{A_j} \right)^{1/q}$. 
\begin{lemma}\label{lem:vector-lq-int} Let  $0<\theta<1$, $1\le q_0,q_1<\infty$ and $q\inv=(1-\theta)q_0\inv+\theta q_1\inv$. 
Let $A_j$ be Banach spaces, $j\in\N$.
  Then \begin{equation}\label{e:ell-q-int}
[\ell^{q_0}(A_j),\ell^{q_1}(B_j)]_\theta
=\ell^{q}([A_j,B_j]_\theta).  \end{equation}
If 
$q_1=\iy$, then 
\begin{equation}\label{e:ell-q-infty-int}
[\ell^{q_0}(A_j),\ell^{\iy}(B_j)]_\theta
=\ell^{q}([A_j,B_j]_\theta)=[\ell^{q_0}(A_j),\ell_0^{\iy}(B_j)]_\theta  ,
\end{equation}
where $\ell_0^{\iy}(B_j):=\{ \{c_j\}\in \ell^{\iy}(B_j): \Vert c_j\Vert_{B_j}\to 0\;as\; j\to \iy\}$.
\end{lemma}

If $1\le q_0,q_1<\iy$, 
(\ref{e:ell-q-int}) also follows from 
Lemma \ref{lem:vector-Lp-int} as a special case where 
the underlying measure space can be taken as $(X,\mu)=\Z$. 
If $
q_1=\iy$, then the remark in \cite[Subsection 1.18.1]{Tr78}  shows that the second statement in 
Lemma \ref{lem:vector-lq-int} is also true.

 In the diagram above 
 in order to  show $S, R$ are continuous mappings,  we  need the following well-known lemma.


\begin{lemma}\label{lem:H-L-h} Let $h(x)$ be a monotonely nonincreasing, radial function 
in $L^1(\R^n)$. Let $h_j(x)=2^{jn/2}h(2^{j/2} x)$ be its scaling.  
Then for all $f$ in $L^1_{loc}(\R^n)$
\[
|\int h_j(x-y)f(y)dy|\le c_n\Vert h\Vert_1 Mf(x),
\]
where $Mf$ denotes the usual Hardy-Littlewood maximal function. 
\end{lemma}
Evidently the decay estimate in (\ref{e:ker-phi-dec}) and Lemma \ref{lem:H-L-h} imply the continuity of 
 $S$ and $R$,  in light of the $L^p(\ell^q)$-valued maximal inequality. 
The proof for $B^{s,q}_p(\cL)$ in part (b)  proceeds in a similar way. 

\subsection{Real interpolation} Peetre's $K$-functional \cite{Pe63} is defined as 
\begin{align*} K(t,a):=K(t,a;A_0,A_1)=\inf (\Vert a_0\Vert_{A_0}+t\Vert a_1\Vert_{A_1}) ,
\end{align*}
where the infimum is taken over all representations of $a=a_0+a_1$,  $a_i\in A_i$. 
Let $0<q\le\iy,0<\theta<1$.  For a given interpolation couple $(A_0,A_1)$, the real interpolation space $(A_0,A_1)_{\theta,q}$
is given by 
\begin{align*} 
(A_0,A_1)_{\theta,q}
=\{a\in A_0+A_1: \Vert a\Vert_{(A_0,A_1)_{\theta,q}}=\left( \int_0^\iy t^{-\theta q} K(t,a)^q \frac{dt}{t}\right)^{1/q}<\iy\}
\end{align*}
with usual modifications if $q=\iy$.


Proof of Theorem \ref{th:F-B-real-int} is similar to 
\cite[Subsection 2.4.2]{Tr83} and \cite[Theorem 6.4.5]{BL76}. 
Define $\ell^{s,q}(A)=\{a=\{a_j\}: \Vert a\Vert_{\ell^{s,q}(A)}=
\Vert \{2^{js}\Vert a_j\Vert_A\}\Vert_{\ell^q}<\iy\}$.
For Besov spaces it follows from  
\begin{align*}
&(\ell^{s_0,q_0}(A_0),\ell^{s_1,q_1}(A_1))_{\theta,q}=\ell^{s,q}((A_0,A_1)_{\theta,q}) ,\label{e:ell-sq-real} 
\end{align*}
 $s=(1-\theta)s_0+\theta s_1 $,
${q}\inv={(1-\theta)}{q_0\inv}+ {\theta}{q_1\inv}$  
 and the commutative diagram for $B_p^{s,q}(\cL)$. 
Consult  \cite{Tr83,Tr78} or
 \cite[Chapter 5,Theorem 6]{Pe76};  both of their proofs  rely on retraction method.
Also see \cite{BL76} for a different proof in the special case involving Sobolev spaces.
In the general case \cite{BL76}  suggests 
using a more concrete characterization of the $K$-functional for the Lorentz space $L^{pq}$.  

For the $F$-space 
  the proof follows from  the commutative diagram for $F_p^{s,q}(\cL)$ and
\begin{align*}
&(L^{p_0}(A_0,w_0),L^{p_1}(A_1,w_1))_{\theta,p}=L^{p}((A_0,A_1)_{\theta,p},w) ,
\end{align*}
where $ {p}\inv={(1-\theta)}{p_0\inv}+ {\theta}{ p_1\inv}$, $w=w_0^{1-\theta}w_1^\theta$,  
 $w_0,w_1$ being two weight functions 
 \cite[Chapter 5
 ]{Pe76}.

\subsection{Schr\"odinger operators with magnetic potential} 
From \cite
{He90a}, \cite{Z06a} 
or \cite{OZ08} 
we know that if the heat kernel of $\cL$ satisfies the upper Gaussian bound 
 \begin{equation}\label{e:g-b-etH}  |e^{-t\cL}(x,y)|\le c_n t^{-n/2} e^{-c |x-y|^2/t}
 \end{equation}
then 
the kernel decay in Assumption \ref{a:phi-dec-Om} holds. 
Let \begin{equation*}
H=  -\sum_{j=1}^n ( \partial_{x_j}+i a_j)^2+ V ,
\end{equation*}  where $a_j(x)\in L^2_{loc}(\R^n)$ is real-valued, 
$V=V_+-V_-$ with $V_+\in L^1_{loc}(\R^n)$, $V_-\in K_n$, the Kato class \cite{Si82}. 
 Proposition 5.1 in 
 \cite{DP05}  
 showed that (\ref{e:g-b-etH}) is valid  for $-\De+V$ 
 if $V_+\in K_n$ and 
$\Vert V_-\Vert_{K_n}< \ga_n:= \pi^{n/2}/\Ga(\frac{n}{2}-1)$, $n\ge 3$, 
whose proof evidently works for $V_+\in L^1_{loc}$. 
By the diamagnetic inequality  \cite[Theorem B.13.2]{Si82}, 
we see that \eqref{e:g-b-etH} also holds for $H$ provided 
$\Vert V_-\Vert_{K_n}< \ga_n$, $n\ge 3$. 

As another example,  a {\em uniformly 
elliptic operator} is given by 
\begin{equation*}
\cL=-\sum_{j,k=1}^n \pa_{x_j}(a_{jk}\pa_{x_k}) ,
\end{equation*}
where $a_{jk}(x)=a_{kj}(x)\in L^\iy(\R^n)$ are real-valued and satisfy the ellipticity condition $(a_{jk})\approx I_n$. 
Then \cite[Theorem 1]{ou06} tells that (\ref{e:g-b-etH}) is true provided that the infimum of its spectrum
$\inf\sigma (\cL)=0$.

\end{document}